\documentclass[epsfig,latexsym,amsfonts,twoside]{article}
\usepackage{amssymb,amsmath,amscd,epsfig}
\usepackage[all,cmtip]{xy}

\def\part#1{\frac{\partial\phantom{#1}}{\partial#1}}
\newtheorem{thm}{Theorem}
\newtheorem{prp}[thm]{Proposition}
\newtheorem{lem}[thm]{Lemma}

\newenvironment{prf}{\begin{trivlist}\item[]{\bf Proof} }%
{\hfill $\Box$ \end{trivlist}}
\newenvironment{dfn}{\begin{trivlist}\item[]{\bf Definition}\em }%
{\end{trivlist}}
\newenvironment{rmk}{\begin{trivlist}\item[]{\bf Remark} }%
{\end{trivlist}}
\newenvironment{exm}{\begin{trivlist}\item[]{\bf Example} }%
{\end{trivlist}}
{\end{trivlist}}

\def\Z{\ifmmode{{\mathbb Z}}\else{${\mathbb Z}$}\fi}
\def\Q{\ifmmode{{\mathbb Q}}\else{${\mathbb Q}$}\fi}
\def\C{\ifmmode{{\mathbb C}}\else{${\mathbb C}$}\fi}
\def\P{\ifmmode{{\mathbb P}}\else{${\mathbb P}$}\fi}

\def\H{\ifmmode{{\mathrm H}}\else{${\mathrm H}$}\fi}

\def\B{\ifmmode{{\cal B}}\else{${\cal B}$}\fi}
\def\E{\ifmmode{{\cal E}}\else{${\cal E}$}\fi}
\def\F{\ifmmode{{\cal F}}\else{${\cal F}$}\fi}
\def\K{\ifmmode{{\cal K}}\else{${\cal K}$}\fi}
\def\L{\ifmmode{{\cal L}}\else{${\cal L}$}\fi}
\def\M{\ifmmode{{\cal M}}\else{${\cal M}$}\fi}
\def\N{\ifmmode{{\cal N}}\else{${\cal N}$}\fi}
\def\O{\ifmmode{{\cal O}}\else{${\cal O}$}\fi}
\def\U{\ifmmode{{\cal U}}\else{${\cal U}$}\fi}
\def\X{\ifmmode{{\cal X}}\else{${\cal X}$}\fi}

\newcommand{\re}{\mathrm{Re}}
\newcommand{\im}{\mathrm{Im}}

\begin{document}

\title{Generalized twistor spaces for hyperk{\"a}hler manifolds}
\author{Rebecca Glover and Justin Sawon}
\date{September, 2013}
\maketitle

\begin{abstract}

Let $M$ be a hyperk{\"a}hler manifold. The $S^2$-family of complex structures compatible with the hyperk{\"a}hler metric can be assembled into a single complex structure on $Z=M\times S^2$; the resulting complex manifold is known as the twistor space of $M$. We describe the analogous construction for generalized complex structures in the sense of Hitchin. Specifically, we exhibit a natural $S^2\times S^2$-family of generalized complex structures compatible with the hyperk{\"a}hler metric, and assemble them into a single generalized complex structure on $\mathcal{X}=M\times S^2\times S^2$. We call the resulting generalized complex manifold the generalized twistor space of $M$.

\end{abstract}

\section{Introduction}

Penrose's twistor theory was developed so that complex analytic techniques could be applied to problems in general relativity, and was later adapted to the Riemannian setting by Atiyah, Hitchin, and Singer~\cite{ahs78}. On each tangent space $T_pM$ of a four-dimensional Riemannian manifold $(M,g)$, the family of complex structures compatible with the metric is parametrized by $\mathrm{SO}(4)/\mathrm{U}(2)\cong S^2$. Globally we get an $S^2$-bundle $Z$ over $M$, and by identifying each fibre with $\mathbb{CP}^1$ we can define a natural almost complex structure on $Z$. Atiyah et al.\ prove that this almost complex structure is integrable if and only if the Riemannian metric on $M$ is self-dual. In this case the complex three-fold $Z$ is known as the {\em twistor space\/} of $M$, and geometric properties of $M$ are completely encoded in the holomorphic structure of $Z$. In particular, the twistor correspondence identifies instantons on $M$ (vector bundles equipped with self-dual connections) with holomorphic bundles on $Z$. This has been one of the most important techniques for constructing and classifying instantons on four-manifolds.

K3 surfaces admit self-dual metrics. In fact, Calabi-Yau metrics on K3 surfaces are hyperk{\"a}hler, which is an even stronger condition. For a hyperk{\"a}hler manifold $(M,g)$ there is an $S^2$-family of global complex structures compatible with the metric, so as a smooth manifold $Z$ is diffeomorphic to $M\times S^2$. However, as a complex three-fold $Z$ is not a product; rather, it is a holomorphic fibration over $\mathbb{CP}^1$. The resulting $\mathbb{CP}^1$-family of complex K3 surfaces is known as a twistor family, and has played an important r\^{o}le in the study of the moduli space of K3 surfaces. Moreover, this construction works just as well in higher dimensions: the geometry of a hyperk{\"a}hler manifold $M$ of real dimension $4n$ is completely encoded in its $(2n+1)$-complex-dimensional twistor space $Z$ (see Hitchin et al.~\cite{hklr87}).

In this paper we extend the twistor construction to generalized complex structures. Hitchin~\cite{hitchin03} defined a generalized complex structure as an orthogonal complex structure on $TM\oplus T^*M$ that satisfies a certain integrability condition. On a hyperk{\"a}hler manifold $M$, we show that there is an $S^2\times S^2$-family of generalized complex structures compatible with the metric. We exhibit this family in three different ways. Firstly, we combine certain natural families, first described by Gualtieri~\cite{gualtieri11}, that interpolate between structures of complex and symplectic type. Secondly, we describe the corresponding family of pure spinors as a quadric in a maximal positive subspace of the Mukai lattice; cf.\ related work of Huybrechts~\cite{huybrechts05}. Thirdly, we use the correspondence between the bi-Hermitian geometry of Gates et al.~\cite{ghr84} and generalized K{\"a}hler geometry introduced by Gualtieri~\cite{gualtieri10}. The reader may prefer to skip to the third approach, which is perhaps the most elegant and revealing. On the other hand, the first approach is really essential for computing a pure spinor (the second approach also yields a pure spinor, and more easily, but only in dimension four).

We then combine these generalized complex structures into a single generalized almost complex structure $\mathcal{J}$ on the smooth manifold $\mathcal{X}=M\times S^2\times S^2$. Our main theorem is that $\mathcal{J}$ is integrable. It is proved by calculating the pure spinor of $\mathcal{J}$ and showing that it is $d$-closed. We call $\mathcal{X}$ (equipped with $\mathcal{J}$) the {\em generalized twistor space\/} of $M$. We also describe some properties of the generalized twistor space $\mathcal{X}$ that are analogues of properties of the usual twistor space $Z$.

Note that the generalized complex structure $\mathcal{J}$ is {\em not\/} of constant type on $\mathcal{X}$. Moreover, it does not seem to arise via previous constructions that yield type changing behaviour, such as blow-ups and Poisson deformations. We hope that a generalized twistor correspondence might eventually be developed, that could lead to new examples of generalized holomorphic bundles.

While working on these results, we discovered a paper of Bredthauer~\cite{bredthauer07} that also gives a construction of generalized twistor spaces for hyperk{\"a}hler manifolds. Bredthauer uses bi-Hermitian geometry, as in our third approach. However, we were unable to follow his proof of the integrability of $\mathcal{J}$: his derivation of the pure spinor is not clearly explained, and the formula itself seems to be incorrect. Pantilie~\cite{pantilie11} constructed twistor spaces for generalized quaternionic manifolds; these come from $\mathbb{CP}^1$-families and sit inside our generalized twistor space. Davidov and Mushkarov~\cite{dm06,dm07} also constructed twistor spaces for generalized complex structures, but their construction is quite different: they do not equip $M$ with a metric, so the fibres of their twistor spaces are always non-compact.

The paper is organized as follows. In Section~2 we recall basic notions from generalized complex geometry. In Section~3 we recall the twistor construction for hyperk{\"a}hler manifolds. In Section~4 we describe our three approaches to families of compatible generalized complex structures on hyperk{\"a}hler manifolds. Then we construct the generalized twistor space and describe its properties.

The authors would like to thank the Hausdorff Research Institute for Mathematics, Bonn, for hospitality during the Junior Trimester Program on differential geometry. The second author gratefully acknowledges support from the Max Planck Institute for Mathematics, Bonn, and from the NSF, grant number DMS-1206309.

\vspace*{3mm}
Conventions: The reader should note that our sign conventions sometimes differ from those in the sources cited, as there is some inconsistency in the literature and it was impossible to achieve internal consistency without altering some signs. Specifically, our parametrization of the $S^2$-family of complex structures on a hyperk{\"a}hler manifold differs from Hitchin et al.~\cite{hklr87}, resulting in a slightly different formula for the holomorphic two-form $\sigma_{\eta}$ in Section~3. In the definition of a generalized K{\"a}hler structure in Section~4.3, our generalized metric $G$ differs from Gualtieri's~\cite{gualtieri10} by a sign, so that the generalized complex structures $(\mathcal{J}_I,\mathcal{J}_{\omega})$ defined in Section~2 yield a genuine example. Consequently, some signs in the formulae of Proposition~\ref{bi-Hermitian} also required alteration.

\section{Generalized complex geometry}

We start with a brief review of generalized complex geometry, a topic initiated by Hitchin~\cite{hitchin03} and further developed by Gualtieri~\cite{gualtieri11}. Let $M$ be a smooth $n$-dimensional manifold with tangent bundle $T$ and cotangent bundle $T^*$. There is a natural inner product on the direct sum $E:=T\oplus T^*$ given by
$$\langle X+\xi,Y+\eta\rangle:=\frac{1}{2}(\xi(Y)+\eta(X)),$$
where $X$ and $Y$ denote tangent vectors and $\xi$ and $\eta$ cotangent vectors. This inner product is indefinite of signature $(n,n)$. There is a natural extension of the Lie bracket of vector fields to smooth sections of $E$, known as the Courant bracket. It is given by
$$[X+\xi,Y+\eta]:=[X,Y]+\mathcal{L}_X\eta-i_Yd\xi$$
where $X$ and $Y$ now denote vector fields and $\xi$ and $\eta$ one-forms. More precisely, this formula gives the Dorfman bracket; the Courant bracket is given by skew-symmetrizing. The (Dorfman) bracket satisfies
\begin{itemize}
\item the Leibniz rule, namely $[e_1,[e_2,e_3]]=[[e_1,e_2],e_3]+[e_2,[e_1,e_3]]$ where $e_i\in C^{\infty}(E)$,
\item $[e_1,e_1]=\pi^*d\langle e_1,e_1\rangle$ where $\pi:E=T\oplus T^*\rightarrow T$ is the projection map.
\end{itemize}
Although it is not skew-symmetric, it becomes skew-symmetric on restricting to isotropic subbundles of $E$, because $\langle e_1,e_1\rangle=0$ implies $[e_1,e_1]=0$.


\begin{rmk}
The Lie bracket of vector fields can be seen as a {\em derived bracket\/} as follows. If $X$ is a vector field, interior product gives a degree $-1$ operator $\iota_X$ on the exterior algebra $C^{\infty}(\wedge^{\bullet}T^*)$. On the other hand, the exterior derivative $d$ gives a degree $1$ operator. The Lie derivative $\mathcal{L}_X$ is then the commutator $[\iota_X,d]$. Now if $Y$ is a second vector field, $[\iota_X,[\iota_Y,d]]$ will be a degree $-1$ operator that looks like $\iota_Z$ for some vector field $Z$. Then $Z=[X,Y]$, and in fact one could take this as the definition of the Lie bracket.

Now $C^{\infty}(E)$ also acts on $C^{\infty}(\wedge^{\bullet}T^*)$, by interior product of vector fields and wedging of one-forms, and the same steps as above lead to a bracket on $C^{\infty}(E)$. This is the bracket described above.
\end{rmk}

The Newlander-Nirenberg Theorem states that an almost complex structure gives a complex structure if and only if it is integrable, meaning that its $+i$-eigenbundle $T^{1,0}\subset T\otimes\mathbb{C}$ is involutive, $[T^{1,0},T^{1,0}]\subset T^{1,0}$. The analogous condition in generalized geometry was used by Hitchin~\cite{hitchin03} to give a definition of generalized complex structures, as follows.

\begin{dfn}
A generalized almost complex structure on $M$ is an orthogonal endomorphism $\mathcal{J}$ of $E=T\oplus T^*$ such that $\mathcal{J}^2$ is minus the identity. Let $L$ denote the $+i$-eigenbundle of $\mathcal{J}$ in $E\otimes\mathbb{C}$. Then $\mathcal{J}$ is a generalized complex structure if $L$ is involutive with respect to the Courant bracket, $[L,L]\subset L$.
\end{dfn}

\begin{rmk}
Note that $L$ is always isotropic, as
$$\langle e,e\rangle=\langle\mathcal{J}e,\mathcal{J}e\rangle=\langle ie,ie\rangle=-\langle e,e\rangle$$
must vanish for $e\in L$. So the Dorfman bracket agrees with its skew-symmetrization, the Courant bracket, when restricted to $L$. If $[L,L]\subset L$ then $L$ becomes a Lie algebroid, known as a {\em (complex) Dirac structure\/}.
\end{rmk}

As mentioned earlier, $C^{\infty}(E)$ acts on $C^{\infty}(\wedge^{\bullet}T^*)$. This action, defined in terms of interior product of vector fields and wedging of one-forms, can actually be defined point-wise. It can also be complexified. Given a generalized almost complex structure $\mathcal{J}$, the canonical line bundle is the complex line subbundle $K\subset\wedge^{\bullet}T^*\otimes\mathbb{C}$ annihilated by the Dirac structure $L\subset E\otimes\mathbb{C}$. A {\em pure spinor\/} for $\mathcal{J}$ is then a (local) generator for $K$, i.e., a (local) section of $\wedge^{\bullet}T^*\otimes\mathbb{C}$ whose annihilator at each point is precisely the fibre of $L$ at that point. Note that $\Phi$ is defined up to multiplication by a smooth complex-valued function, and $\mathcal{J}$ is completely determined (locally) by $\Phi$. Integrability of $\mathcal{J}$ was formulated in terms of pure spinors by Gualtieri, Theorem~2.9 of~\cite{gualtieri11}.

\begin{prp}
\label{integrable}
A generalized almost complex structure $\mathcal{J}$ is integrable if and only if it is represented at every point by a pure spinor $\Phi$ that satisfies
$$d\Phi=\iota_X\Phi+\xi\wedge\Phi$$
for some (local) section $X+\xi\in C^{\infty}(E\otimes\mathbb{C})$. In particular, the vanishing of $d\Phi$ is a sufficient condition for integrability.
\end{prp}

\begin{exm}
If $M$ admits an almost complex structure $I:T\rightarrow T$ then
$$\mathcal{J}_I:=\left(\begin{array}{cc} -I & 0 \\ 0 & I^* \end{array}\right)$$
is a generalized almost complex structure on $M$. One can show that integrability of $\mathcal{J}_I$ is equivalent to integrability of $I$. The Dirac structure for $\mathcal{J}_I$ is $T^{0,1}\oplus\Omega^{1,0}$, and pure spinors are given by $e_1\wedge\ldots\wedge e_n\in C^{\infty}(\wedge^nT^*\otimes\mathbb{C})$ where $\{e_1,\ldots,e_n\}$ are local bases for $(1,0)$-forms. For example, if $(z_1,\ldots,z_n)$ are local holomorphic coordinates on $M$ then $dz_1\wedge\ldots\wedge dz_n$ is a pure spinor.
\end{exm}

\begin{exm}
If $M$ admits an almost symplectic structure $\omega$ (a non-degenerate two-form, not necessarily closed) then
$$\mathcal{J}_{\omega}:=\left(\begin{array}{cc} 0 & -\omega^{-1} \\ \omega & 0 \end{array}\right)$$
is a generalized almost complex structure on $M$, where we think of $\omega$ as a map from $T$ to $T^*$ given by interior product. One can show that $\mathcal{J}_{\omega}$ is integrable if and only if $\omega$ is a symplectic structure, i.e., $d\omega=0$. The Dirac structure for $\mathcal{J}_{\omega}$ is
$$\{X-i\omega(X)\,|\,X\in C^{\infty}(T)\},$$
and a pure spinor is given by $\mathrm{exp}(i\omega)\in C^{\infty}(\wedge^{\mathrm{even}}T^*\otimes\mathbb{C})$.
\end{exm}

The Lie bracket of vector fields is invariant under diffeomorphisms, and in fact the automorphisms of the tangent bundle $T$ equipped with the Lie bracket on its sections are precisely the diffeomorphisms of $M$. On the other hand, $T\oplus T^*$ equipped with the Courant bracket on its sections has an enlarged automorphism group: the additional symmetries come from $B$-field transforms.

\begin{dfn}
Let $B$ be a closed two-form, regarded as a map from $T$ to $T^*$ given by interior product. Then
$$e^B:=\left(\begin{array}{cc} 1 & 0 \\ B & 1 \end{array}\right):T\oplus T^*\rightarrow T\oplus T^*$$
is a symmetry of $T\oplus T^*$ that preserves the Courant bracket, known as a $B$-field transform.
\end{dfn}

These symmetries can be used to generate new generalized complex structures, as we demonstrate below.

\begin{exm}
Suppose that $M$ admits a generalized complex structure $\mathcal{J}$ with Dirac structure $L$ and pure spinor $\Phi$. Then the $B$-field transform $\mathcal{J}_B:=e^{-B}\mathcal{J}e^B$ is also a generalized complex structure, with Dirac structure $e^{-B}L$ and pure spinor $e^B\Phi$. In particular, when $\mathcal{J}=\mathcal{J}_{\omega}$ comes from a symplectic structure on $M$, a pure spinor for the $B$-field transform $\mathcal{J}_B$ is given by $\mathrm{exp}(B+i\omega)$.
\end{exm}

\begin{dfn}
A generalized complex manifold of real dimension $2n$ (the dimension must be even) is stratified by type, which at the point $p\in M$ is given by
$$\frac{1}{2}\mathrm{dim}(T^*_pM\cap\mathcal{J}T^*_pM).$$
\end{dfn}

The type gives an upper semi-continuous function on $M$; it takes values in $\{0,1,\ldots,n\}$ and can only change by an even integer (we assume $M$ is connected). If the type is constant and equal to $n$ at every point then $\mathcal{J}$ must come from a complex structure $I$, and we say that $\mathcal{J}=\mathcal{J}_I$ is of {\em complex type\/}. If the type is constant and equal to $0$ at every point then $\mathcal{J}$ must either come from a symplectic structure $\omega$ or the $B$-field transform of a symplectic structure; in either case, we say $\mathcal{J}=\mathcal{J}_{\omega}$ or $e^{-B}\mathcal{J}_{\omega}e^B$ is of {\em symplectic type\/}.

The ``Generalized Darboux Theorem'' of Gualtieri~\cite{gualtieri11} gives a local classification of generalized complex manifolds of constant type: on an open set where the type is constant and equal to $k$, a generalized complex manifold is equivalent to the product of an open set in $\mathbb{C}^k$ and an open set in $\mathbb{R}^{2n-2k}$ with its standard symplectic structure.

More interesting generalized complex manifolds exhibit type-jumping. Some methods for constructing such examples include blowing-up, reduction, Poisson deformation, and the bi-Hermitian/generalized K{\"a}hler correspondence. Bailey~\cite{bailey12} gave a complete local classification by proving that a generalized complex manifold is locally equivalent to the product of a symplectic manifold and a Poisson deformation of a complex manifold. In this article we give a new construction of generalized complex manifolds exhibiting type-jumping. The following example of Gualtieri~\cite{gualtieri11} is fundamental to our construction.

\begin{exm}
Let $M$ be a holomorphic symplectic manifold with complex structure $I$ and holomorphic symplectic form $\sigma$. Write $\sigma=\omega_J+i\omega_K$ for the decomposition into real and imaginary parts. It is easy to check that
$$\mathcal{J}_{\theta}:=\cos\theta\left(\begin{array}{cc} -I & 0 \\ 0 & I^* \end{array}\right)+\sin\theta\left(\begin{array}{cc} 0 & -\omega_J^{-1} \\ \omega_J & 0 \end{array}\right)$$
gives a generalized almost complex structure for all $\theta\in\mathbb{R}$. In fact, $\mathcal{J}_{\theta}$ is always integrable. This is obvious when $\theta$ is a multiple of $\pi$, as $\mathcal{J}_{\theta}$ is then $\mathcal{J}_I$ or $\mathcal{J}_{-I}$. For all other $\theta$,
$$\mathcal{J}_{\theta}=e^{-B}\mathcal{J}_{\omega}e^B$$
is the $B$-field transform of a symplectic structure, where $B=-(\cot\theta)\omega_K$ and $\omega=(\csc\theta)\omega_J$.

As described by the second author in~\cite{sawon12}, this example can be extended to a $\mathbb{CP}^1$-family of generalized complex structures. We first replace $\omega_J$ by a combination $(\cos\phi)\omega_J+(\sin\phi)\omega_K$, which gives
$$\mathcal{J}_{\theta,\phi}:=(\cos\theta)\mathcal{J}_I+(\sin\theta\cos\phi)\mathcal{J}_{\omega_J}+(\sin\theta\sin\phi)\mathcal{J}_{\omega_K}.$$
Regarding $\theta$ and $\phi$ as spherical coordinates on the sphere $S^2$, and then using stereographic projection to change to a single complex parameter $\zeta$ on the extended complex plane $\mathbb{CP}^1$ gives
$$\mathcal{J}_{\zeta}:=\left(\frac{1-|\zeta|^2}{1+|\zeta|^2}\right)\mathcal{J}_I+\left(\frac{2\re\zeta}{1+|\zeta|^2}\right)\mathcal{J}_{\omega_J}+\left(\frac{2\im\zeta}{1+|\zeta|^2}\right)\mathcal{J}_{\omega_K}.$$
The pure spinor defining $\mathcal{J}_{\zeta}$ is given by
$$\Phi_{\zeta}=(2i\zeta)^n\exp\left(\frac{\sigma}{2i\zeta}-\frac{i\zeta\bar{\sigma}}{2}\right)$$
where $M$ has real dimension $4n$. The factor of $(2i\zeta)^n$ is included to clear denominators (for $k>n$, note that there are no $(2k,0)$-forms, and hence $\sigma^k=0$) and give a leading term of $\sigma^n$ when $\zeta=0$. For instance, when $n=1$ expanding gives
$$\Phi_{\zeta}=\sigma+2i\zeta\left(1-\frac{1}{4}\sigma\bar{\sigma}\right)+\zeta^2\bar{\sigma}.$$
In general, a crucial feature of the formula for $\Phi_{\zeta}$ is that it depends holomorphically on $\zeta$, so that we obtain a holomorphic family of generalized complex structures on $M$.
\end{exm}

On a hyperk{\"a}hler manifold we have a triple of complex structures $I$, $J$, and $K$, and a triple of corresponding K{\"a}hler forms $\omega_I$, $\omega_J$, and $\omega_K$. Our first goal is to include the resulting six generalized complex structures $\mathcal{J}_I$, $\mathcal{J}_J$, $\mathcal{J}_K$, $\mathcal{J}_{\omega_I}$, $\mathcal{J}_{\omega_J}$, and $\mathcal{J}_{\omega_K}$ in a single holomorphic family. Before that, we review the usual twistor construction which involves only complex structures.

\section{Twistor spaces for hyperk{\"a}hler manifolds}

\begin{dfn}
A hyperk{\"a}hler manifold is a Riemannian manifold $M$ whose metric $g$ is K{\"a}hlerian with respect to a triple of complex structures $I$, $J$, and $K$ that behave like the quaternions, i.e.,
$$IJ=K=-JI,\qquad JK=I=-KJ,\qquad\mbox{and}\qquad KI=J=-IK.$$
We denote the corresponding K{\"a}hler forms by $\omega_I$, $\omega_J$, and $\omega_K$. The real dimension of $M$ is necessarily a multiple of four, $\mathrm{dim}_{\mathbb{R}}M=4n$.
\end{dfn}

It is easy to check that $g$ is also K{\"a}hlerian with respect to the complex structure $aI+bJ+cK$, where $a^2+b^2+c^2=1$, i.e., $(a,b,c)\in S^2$. In fact, this is the complete family of complex structures compatible with the hyperk{\"a}hler metric. Identifying $S^2$ with $\mathbb{CP}^1$ and changing to a single complex parameter $\eta$ gives
$$I_{\eta}=\left(\frac{1-|\eta|^2}{1+|\eta|^2}\right)I+\left(\frac{2\re\eta}{1+|\eta|^2}\right)J+\left(\frac{2\im\eta}{1+|\eta|^2}\right)K.$$
A corresponding holomorphic two-form for $I_{\eta}$ is given by
$$\sigma_{\eta}=(\omega_J+i\omega_K)-2\eta\omega_I-\eta^2(\omega_J-i\omega_K).$$
Of course, for each $\eta\in\mathbb{CP}^1$, the holomorphic two-form is only defined up to multiplication by a non-zero factor; the above normalization is chosen so that $\sigma_{\eta}$ depends holomorphically on $\eta$, which will again be crucial.

The following lemma, proved by a direct calculation, will also be useful later.

\begin{lem}
\label{SO(3)}
Write $\eta=x+iy$, let $\sigma_{\eta}^{\prime}:=\sigma_{\eta}/(1+x^2+y^2)$, and let $\omega_{\eta}$ denote the K{\"a}hler form of $I_{\eta}$. Then
$$\left(\begin{array}{c} \omega_{\eta} \\ \re\sigma_{\eta}^{\prime} \\ \im\sigma_{\eta}^{\prime} \end{array}\right)=\frac{1}{1+x^2+y^2}
\left(\begin{array}{ccc} 1-x^2-y^2 & 2x & 2y \\
-2x & 1-x^2+y^2 & -2xy \\
-2y & -2xy & 1+x^2-y^2 \end{array}\right)\left(\begin{array}{c} \omega_I \\ \omega_J \\ \omega_K \end{array}\right)$$
where, moreover, the coefficient matrix lies in $\mathrm{SO}(3)$.
\end{lem}

With respect to $I_{\eta}$, $\sigma_{\eta}^n$ is a $(2n,0)$-form. Thus the $-i$-eigenspace of $I_{\eta}$,
$$T^{0,1}_{\eta}\subset T\otimes\mathbb{C},$$
consists of vectors whose interior product with $\sigma_{\eta}^n$ is zero. In this way, $\sigma_{\eta}$ determines the complex structure $I_{\eta}$; indeed, in the language of generalized complex structures, $\sigma_{\eta}^n$ is a pure spinor for $\mathcal{J}_{I_{\eta}}$.

We can define an almost complex structure on the smooth manifold $M\times\mathbb{CP}^1$ by
$$(I_{\eta},I_{\mathbb{CP}^1})\in\mathrm{End}\left(T_{(m,\eta)}M\times\mathbb{CP}^1\right)=\mathrm{End}\left(T_mM\right)\times \mathrm{End}\left(T_{\eta}\mathbb{CP}^1\right),$$
where $I_{\mathbb{CP}^1}$ denotes the standard complex structure on $\mathbb{CP}^1$. The pure spinor for this almost complex structure is $\sigma_{\eta}^n\wedge d\eta$. Since
$$d(\sigma_{\eta}^n\wedge d\eta)=n\sigma_{\eta}^{(n-1)}\wedge\left(\frac{\partial\sigma_{\eta}}{\partial\eta}\right)d\eta\wedge d\eta=0,$$
the almost complex structure is integrable.

\begin{dfn}
The twistor space of $M$ is the complex manifold $Z$ of dimension $2n+1$ given by putting the complex structure $(I_{\eta},I_{\mathbb{CP}^1})$ on the smooth manifold $M\times\mathbb{CP}^1$.
\end{dfn}

The twistor space comes equipped with various structures which together encode the hyperk{\"a}hler metric on $M$; see Theorem~3.3 of Hitchin et al.~\cite{hklr87}.

\begin{thm}
\label{twistor}
The twistor space $Z$ has the following properties:
\begin{enumerate}
\item it admits a holomorphic fibration $p:Z\rightarrow\mathbb{CP}^1$,
\item $p$ admits a family of holomorphic sections, each with normal bundle $\mathcal{O}(1)^{\oplus 2n}$,
\item there is a holomorphic section of $\Lambda^2T^*_F(2)$, where $T_F$ denotes the tangent bundle to the fibres, defining a holomorphic symplectic form on each fibre of $p$,
\item $Z$ has a real structure compatible with the other structures and inducing the antipodal map on $\mathbb{CP}^1$.
\end{enumerate}
Conversely, given a complex manifold $Z$ of complex dimension $2n+1$ with the above structures, the space parametrizing real twistor lines, i.e., real sections of $p:Z\rightarrow\mathbb{CP}^1$, is a manifold of real dimension $4n$ with a natural hyperk{\"a}hler metric whose twistor space is $Z$.
\end{thm}


\section{Generalized twistor spaces}

Let $(M,g)$ be a hyperk{\"a}hler manifold. The metric defines a generalized metric on $E=T\oplus T^*$,
$$G(X+\xi,Y+\eta):=\frac{1}{2}\left(g(X,Y)+g^{-1}(\xi,\eta)\right).$$
The generalized complex structures $\mathcal{J}$ compatible with $G$, in the sense that $G(\mathcal{J}e_1,\mathcal{J}e_2)=G(e_1,e_2)$ for $e_i\in C^{\infty}(E)$, will be parametrized by $\mathbb{CP}^1\times\mathbb{CP}^1$. We begin this section by describing this family of generalized complex structures in three different ways.

\subsection{Hyperk{\"a}hler rotation}

Recall the example at the end of Section~2: for a holomorphic symplectic manifold $M$ with complex structure $I$ and holomorphic symplectic form $\sigma$ we constructed a $\mathbb{CP}^1$-family of generalized complex structures $\mathcal{J}_{\zeta}$. We can apply this construction to a hyperk{\"a}hler manifold, because $\sigma=\omega_J+i\omega_K$ is a holomorphic symplectic form with respect to $I$. More generally,
$$\sigma_{\eta}=(\omega_J+i\omega_K)-2\eta\omega_I-\eta^2(\omega_J-i\omega_K)$$
is a holomorphic symplectic form with repect to $I_{\eta}$. So for each $\eta\in\mathbb{CP}^1$ we obtain a family $\mathcal{J}_{\eta,\zeta}$ of generalized complex structures on $M$ parametrized by $\zeta\in\mathbb{CP}^1$. Stated equivalently, the family $\mathcal{J}_{\eta,\zeta}$ is parametrized by a $\mathbb{CP}^1$-bundle over $\mathbb{CP}^1$. We can identify this bundle.

\begin{lem}
Associated to each hyperk{\"a}hler metric $g$, there is a family of generalized complex structures on $M$ parametrized by the projectivization
$$\mathbb{P}(T\oplus\O)\longrightarrow\mathbb{CP}^1$$
of the tangent bundle on $\mathbb{CP}^1$.
\end{lem}

\begin{prf}
Following the above discussion, it remains to identify the $\mathbb{CP}^1$-bundle. The important point is that the holomorphic symplectic form $\sigma_{\eta}$ is not well-defined for all $\eta\in\mathbb{CP}^1$; it has a pole at $\infty$. The correct way to proceed is to think of $\sigma_{\eta}$ as a section of the bundle $\Lambda^2T^*_{\mathrm{fibre}}(2)$ over $\mathbb{CP}^1$, as in Theorem~\ref{twistor}. Then locally on the base $\mathbb{CP}^1$, $\sigma_{\eta}$ can be multiplied by a non-vanishing section of $\mathcal{O}(-2)$ to produce a family of holomorphic symplectic forms, each generating a $\mathbb{CP}^1$-family of generalized complex structures as above. Therefore the $\mathbb{CP}^1$-bundle is the projectivization of $\mathcal{O}(-2)$. This is isomorphic to
$$\mathbb{P}(\mathcal{O}\oplus\mathcal{O}(-2))\cong\mathbb{P}(\mathcal{O}(2)\oplus\mathcal{O})\cong\mathbb{P}(T\oplus\mathcal{O}).$$
\end{prf}

Now, there is duplication in this family. Firstly, let $\eta$ and $-\bar{\eta}^{-1}$ be antipodal points in the base $\mathbb{CP}^1$. Since $I_{-\bar{\eta}^{-1}}=-I_{\eta}$, the $\mathbb{CP}^1$-fibres above these two points yield the same family of generalized complex structures: in one fibre, the north and south poles correspond to $\mathcal{J}_{I_{\eta}}$ and $\mathcal{J}_{-I_{\eta}}$, while in the other fibre these are reversed. Secondly, there is additional overlap of the equators of the fibres. Recall that the equator of the fibre above $0$ parametrizes generalized complex structures of the form
$$\mathcal{J}_{(\cos\phi)\omega_J+(\sin\phi)\omega_K}.$$
A generic equator will also parametrize generalized complex structures coming from symplectic structures, i.e., of the form
$$\mathcal{J}_{a\omega_I+b\omega_J+c\omega_K}$$
where $a^2+b^2+c^2=1$. However, there is only an $S^2$-family of such symplectic structures, whereas the equators give an $S^1$-bundle over $\mathbb{CP}^1$ (this is actually the Hopf fibration on $S^3$).

We will show that this family is actually pulled back from a $\mathbb{CP}^1\times\mathbb{CP}^1$-family via a map that we now describe. Although we will describe this map locally using `complex' coordinates, it is not a holomorphic map, as it will be derived from the real exponential map.

\begin{lem}
\label{map_f}
Define a local coordinate patch on $\mathbb{P}(T\oplus\mathcal{O})$ by
\begin{eqnarray*}
\mathbb{C}^2 & \longrightarrow & T\mathbb{C}\subset T\mathbb{CP}^1\subset\mathbb{P}(T\oplus\mathcal{O}) \\
(\eta,\zeta) & \longmapsto & \zeta(1+|\eta|^2)\frac{\partial\phantom{\eta}}{\partial\eta}.
\end{eqnarray*}
Then there exists a map
$$f:\mathbb{P}(T\oplus\mathcal{O})\longrightarrow\mathbb{CP}^1\times\mathbb{CP}^1$$
that is given in these local coordinates by
$$(\eta,\zeta)\longmapsto\left(\frac{\zeta+\eta}{-\bar{\eta}\zeta+1},\frac{-\zeta+\eta}{\bar{\eta}\zeta+1}\right),$$
where on the right hand side we think of the coordinates as lying in $\mathbb{C}\cup\{\infty\}=\mathbb{CP}^1$.
\end{lem}

\begin{rmk}
Note that the local `complex' coordinates $(\eta,\zeta)$ on $T\mathbb{C}$ are not holomorphic, since $(1+|\eta|^2)\frac{\partial\phantom{\eta}}{\partial\eta}$ gives a smooth trivialization of the complex bundle $T\mathbb{C}$, but it does not give a holomorphic trivialization of the holomorphic bundle $T\mathbb{C}$. On the other hand, these coordinates arise naturally when we replace $\sigma_{\eta}$ by the renormalized form $\sigma^{\prime}_{\eta}=\sigma_{\eta}/(1+|\eta|^2)$ of Lemma~\ref{SO(3)}.
\end{rmk}

\begin{prf}
We start with a differential geometric map, identifying $\mathbb{CP}^1$ with the unit sphere $S^2$ with the standard metric. The normal bundle $N_{\Delta}$ to the diagonal $\Delta$ in $S^2\times S^2$ can be canonically identified with its tangent bundle, by
\begin{eqnarray*}
T\Delta & \stackrel{\cong}{\longrightarrow} & N_{\Delta}\subset TS^2\times TS^2 \\
v & \longmapsto & (v,-v).
\end{eqnarray*}
Giving $S^2\times S^2$ the product metric, we then have the exponential map
\begin{eqnarray*}
T\Delta & \longrightarrow & S^2\times S^2 \\
v\in T_p\Delta & \longmapsto & \exp_{(p,p)}(v,-v)=\gamma(1),
\end{eqnarray*}
where $\gamma$ is a geodesic in $S^2\times S^2$ with initial point $\gamma(0)=(p,p)$ and initial direction $\gamma^{\prime}(0)=(v,-v)$. For fixed $p\in\Delta$, the exponential map will be injective up to radius $\pi$, then it will map all vectors $v$ of length $\pi$ to the same point $(-p,-p)\in S^2\times S^2$, where $-p$ denotes the point antipodal to $p$. Let $B\Delta\subset T\Delta$ denote the bundle consisting of open discs of radius $\pi$, let $\bar{B}\Delta\subset T\Delta$ denote its closure, and let $B^*\Delta$ denote the $S^2$-bundle obtained from $\bar{B}\Delta$ by identifying the vectors of length $\pi$ in each fibre. Then the exponential map induces a map
$$B^*\Delta\longrightarrow S^2\times S^2.$$
This map is not injective, though it is injective on each fibre. Finally, we identify $\Delta$ with $\mathbb{CP}^1$, we identify $B\Delta$ with $T\mathbb{CP}^1$ via stereographic projection (note that for each $p$, $B_p\Delta$ is the two-sphere $B^*_p\Delta$ minus the `north pole'), and we identify $B^*\Delta$ with $\mathbb{P}(T\oplus\mathcal{O})$ by compactifying each fibre. This gives the map
$$\mathbb{P}(T\oplus\mathcal{O})\longrightarrow S^2\times S^2=\mathbb{CP}^1\times\mathbb{CP}^1$$
that we call $f$. Note that on the right hand side, the identification $S^2=\mathbb{CP}^1$ also comes from stereographic projection onto $\mathbb{C}\subset\mathbb{CP}^1$, followed by compactification.

Now introduce the local coordinates $(\eta,\zeta)$ on $T\mathbb{C}\subset T\mathbb{CP}^1\subset\mathbb{P}(T\oplus\mathcal{O})$, representing the tangent vector $\zeta(1+|\eta|^2)\frac{\partial\phantom{\eta}}{\partial\eta}$. The point $\eta=0\in\mathbb{CP}^1$ corresponds to $(1,0,0)\in\Delta=S^2$. At this point, the exponential map
$$B^*_{(1,0,0)}\Delta\longrightarrow S^2\times S^2$$
takes a point $(x,y,z)$ in the sphere $B^*_{(1,0,0)}\Delta\subset\mathbb{R}^3$ to $((x,y,z),(x,-y,-z))$ in $S^2\times S^2$. Conjugating by stereographic projection shows that
$$f|_{T_0\mathbb{CP}^1}:T_0\mathbb{CP}^1\longrightarrow B_{(1,0,0)}\Delta\longrightarrow S^2\times S^2\longrightarrow\mathbb{CP}^1\times\mathbb{CP}^1$$
takes $\zeta\frac{\partial\phantom{\eta}}{\partial\eta}$ to $(\zeta,-\zeta)$.

To determine the map $f$ for general $\eta$, we use the fact that $\mathrm{PSU}(2)\cong\mathrm{SO}(3)$ acts transitively on $\mathbb{CP}^1$, and the map $f$ is equivariant by construction. There is a unique group element
$$A:=\frac{1}{\sqrt{1+|\eta|^2}}\left(\begin{array}{cc} 1 & \eta \\ -\bar{\eta} & 1 \end{array}\right)\in\mathrm{PSU}(2)$$
that takes $0$ to $\eta$, and such that $DA_0$ takes $\frac{\partial\phantom{\eta}}{\partial\eta}\in T_0\mathbb{CP}^1$ to $(1+|\eta|^2)\frac{\partial\phantom{\eta}}{\partial\eta}\in T_{\eta}\mathbb{CP}^1$. Applying this map to $f|_{T_0\mathbb{CP}^1}$ shows that
$$f|_{T_{\eta}\mathbb{CP}^1}:T_{\eta}\mathbb{CP}^1\longrightarrow\mathbb{CP}^1\times\mathbb{CP}^1$$
takes $\zeta (1+|\eta|^2)\frac{\partial\phantom{\eta}}{\partial\eta}$ to
\begin{eqnarray*}
(A.\zeta,A.(-\zeta)) & = & \left(A\left(\begin{array}{c} \zeta \\ 1 \end{array}\right),A\left(\begin{array}{c} -\zeta \\ 1 \end{array}\right)\right) \\
 & = & \left(\frac{1}{\sqrt{1+|\eta|^2}}\left(\begin{array}{c} \zeta+\eta \\ -\bar{\eta}\zeta+1 \end{array}\right),\frac{1}{\sqrt{1+|\eta|^2}}\left(\begin{array}{c} -\zeta+\eta \\ \bar{\eta}\zeta+1 \end{array}\right)\right) \\
 & = & \left(\frac{\zeta+\eta}{-\bar{\eta}\zeta+1},\frac{-\zeta+\eta}{\bar{\eta}\zeta+1}\right).
\end{eqnarray*}
This completes the proof.
\end{prf}

\begin{rmk}
Inside $\mathrm{SO}(3)$, the isotropy subgroup of any point in $S^2$ is isomorphic to $S^1$, and therefore there is an $S^1$-family of group elements taking $0$ to $\eta$ in $\mathbb{CP}^1=S^2$. By adding tangent vectors, we can produce orthonormal frames in $\mathbb{R}^3$,
$$e_1=0\in S^2,\qquad e_2=\frac{\partial\phantom{\eta}}{\partial\eta}\in T_0S^2\cong e_1^{\perp},\qquad\mbox{and}\qquad e_3=e_1\times e_2$$
and
$$e^{\prime}_1=\eta\in S^2,\qquad e^{\prime}_2=(1+|\eta|^2)\frac{\partial\phantom{\eta}}{\partial\eta}\in T_\eta S^2\cong (e^{\prime}_1)^{\perp},\qquad\mbox{and}\qquad e^{\prime}_3=e^{\prime}_1\times e^{\prime}_2.$$
Here $(1+|\eta|^2)$ is a normalization factor to make $e^{\prime}_2$ a unit vector. Because $\mathrm{PSU}(2)\cong\mathrm{SO}(3)$ acts faithfully on orthonormal frames in $\mathbb{R}^3$, there is a unique $A\in\mathrm{PSU}(2)$ taking $\{e_1,e_2,e_3\}$ to $\{e^{\prime}_1,e^{\prime}_2,e^{\prime}_3\}$, namely the matrix
$$A=\frac{1}{\sqrt{1+|\eta|^2}}\left(\begin{array}{cc} 1 & \eta \\ -\bar{\eta} & 1 \end{array}\right)$$
appearing in the proof. Note that under the isomorphism $\mathrm{PSU}(2)\cong\mathrm{SO}(3)$, $A$ corresponds to the matrix
$$\frac{1}{1+x^2+y^2}
\left(\begin{array}{ccc} 1-x^2-y^2 & 2x & 2y \\
-2x & 1-x^2+y^2 & -2xy \\
-2y & -2xy & 1+x^2-y^2 \end{array}\right)$$
from Lemma~\ref{SO(3)}, where $\eta=x+iy$.
\end{rmk}

Let us return now to our family $\mathcal{J}_{\eta,\zeta}$ of generalized complex structures on $M$, parametrized by $\mathbb{P}(T\oplus\mathcal{O})$. We determine a pure spinor for this family.

\begin{lem}
A pure spinor for the generalized complex structure $\mathcal{J}_{\eta,\zeta}$ is
\begin{eqnarray*}
\Phi_{\eta,\zeta} & = & (2i\zeta)^n\exp\left(\frac{\sigma^{\prime}_{\eta}}{2i\zeta}-\frac{i\zeta\bar{\sigma}^{\prime}_{\eta}}{2}\right) \\
 & = & (2i\zeta)^n\exp\left(\frac{\sigma-2\eta\omega_I-\eta^2\bar{\sigma}}{2i\zeta(1+|\eta|^2)}-\frac{i\zeta(\bar{\sigma}-2\bar{\eta}\omega_I-\bar{\eta}^2\sigma)}{2(1+|\eta|^2)}\right).
\end{eqnarray*}
\end{lem}

\begin{prf}
This is just the same formula as before, namely
$$\Phi_{\zeta}=(2i\zeta)^n\exp\left(\frac{\sigma}{2i\zeta}-\frac{i\zeta\bar{\sigma}}{2}\right),$$
except now we let $\sigma$ vary too, by replacing it with
$$\sigma^{\prime}_{\eta}=\frac{\sigma-2\eta\omega_I-\eta^2\bar{\sigma}}{1+|\eta|^2}.$$
Recall that the renormalized form $\sigma^{\prime}_{\eta}$ is the appropriate form to use with the coordinates $(\eta,\zeta)$, rather than $\sigma_{\eta}$, as mentioned earlier in the remark following Lemma~\ref{map_f}.
\end{prf}

Using this formula, we can show that the family $\mathcal{J}_{\eta,\zeta}$ is pulled back from $\mathbb{CP}^1\times\mathbb{CP}^1$.

\begin{prp}
Write
$$(\alpha,\beta)=\left(\frac{\zeta+\eta}{-\bar{\eta}\zeta+1},\frac{-\zeta+\eta}{\bar{\eta}\zeta+1}\right)$$
for the image of $(\eta,\zeta)$ under $f$. Then the pure spinor $\Phi_{\eta,\zeta}$ can be written in terms of $\alpha$ and $\beta$ (at least after rescaling), and therefore the family $\mathcal{J}_{\eta,\zeta}$ is the pull back by $f$ of a family $\mathcal{J}_{\alpha,\beta}$ of generalized complex structures on $M$ parametrized by $\mathbb{CP}^1\times\mathbb{CP}^1$.
\end{prp}

\begin{prf}
We begin by computing
$$\alpha\beta=\frac{\eta^2-\zeta^2}{1-\bar{\eta}^2\zeta^2},\qquad \alpha+\beta=\frac{2\eta+2\bar{\eta}\zeta^2}{1-\bar{\eta}^2\zeta^2},\qquad\mbox{and}\qquad \alpha-\beta=\frac{2\zeta(1+|\eta|^2)}{1-\bar{\eta}^2\zeta}.$$
Therefore
\begin{eqnarray*}
\frac{\sigma-(\alpha+\beta)\omega_I-\alpha\beta\bar{\sigma}}{i(\alpha-\beta)} & = & \left(\frac{1-\bar{\eta}^2\zeta^2}{2i\zeta(1+|\eta|^2)}\right)\sigma+\left(\frac{-2\eta-2\bar{\eta}\zeta^2}{2i\zeta(1+|\eta|^2)}\right)\omega_I+\left(\frac{-\eta^2+\zeta^2}{2i\zeta(1+|\eta|^2)}\right)\bar{\sigma} \\
 & = & \frac{\sigma-2\eta\omega_I-\eta^2\bar{\sigma}}{2i\zeta(1+|\eta|^2)}-\frac{i\zeta(\bar{\sigma}-2\bar{\eta}\omega_I-\bar{\eta}^2\sigma)}{2(1+|\eta|^2)} \\
 & = & \frac{\sigma^{\prime}_{\eta}}{2i\zeta}-\frac{i\zeta\bar{\sigma}^{\prime}_{\eta}}{2}.
\end{eqnarray*}
So up to a factor, we see that $\Phi_{\eta,\zeta}$ is equal to
$$\exp\left(\frac{\sigma-(\alpha+\beta)\omega_I-\alpha\beta\bar{\sigma}}{i(\alpha-\beta)}\right).$$
Of course, since the generalized complex structure $\mathcal{J}_{\eta,\zeta}$ only determines a pure spinor line, we are free to rescale $\Phi_{\eta,\zeta}$ as needed.
\end{prf}

\begin{dfn}
Define a pure spinor
$$\Phi_{\alpha,\beta}:=i^n(\alpha-\beta)^n\exp{\left(\frac{\sigma-(\alpha+\beta)\omega_I-\alpha\beta\bar{\sigma}}{i(\alpha-\beta)}\right)},$$
and let the corresponding $\mathbb{CP}^1\times\mathbb{CP}^1$-family of generalized complex structures on $M$ be denoted $\mathcal{J}_{\alpha,\beta}$. (Then the above proof shows that $\mathcal{J}_{\eta,\zeta}$ is the pull back of $\mathcal{J}_{\alpha,\beta}$.)
\end{dfn}

\begin{rmk}
The fact that $\Phi_{\alpha,\beta}$ is meromorphic in $\alpha$ and $\beta$ will be of fundamental importance. In fact, it is holomorphic, as we now show.
\end{rmk}

Expanding the exponential, we get
$$\Phi_{\alpha,\beta}=\sum_{j=0}^{2n}\frac{1}{j!}i^{n-j}(\alpha-\beta)^{n-j}(\sigma-(\alpha+\beta)\omega_I-\alpha\beta\bar{\sigma})^j.$$
A priori, the terms for $j>n$ contain poles along the diagonal $\alpha=\beta$.

\begin{lem}
\label{alpha=beta}
If $j\neq n$ then the $j$th term in $\Phi_{\alpha,\beta}$ is divisible by $(\alpha-\beta)$. Therefore, when $\alpha=\beta$ we have
$$\Phi_{\alpha,\alpha}=\frac{1}{n!}(\sigma-2\alpha\omega_I-\alpha^2\bar{\sigma})^n.$$
In other words, along the diagonal $\Phi_{\alpha,\alpha}$ is indeed the correct pure spinor for the generalized complex structure $\mathcal{J}_{\alpha,\alpha}=\mathcal{J}_{I_{\alpha}}$ of complex type.
\end{lem}

\begin{prf}
The statement is trivial for $j<n$, so assume $j>n$, and write $j=n+k$ with $0<k\leq n$. Denote $\sigma-(\alpha+\beta)\omega_I-\alpha\beta\bar{\sigma}$ by $\sigma_{\alpha,\beta}$ and $\sigma-2\alpha\omega_I-\alpha^2\bar{\sigma}$ by $\sigma_{\alpha}$. Then
$$\sigma_{\alpha,\beta}=\sigma_{\alpha}+(\alpha-\beta)\tau_{\alpha}$$
where $\tau_{\alpha}=\omega_I+\alpha\bar{\sigma}$ is $\frac{\partial\phantom{\alpha}}{\partial\alpha}$ of $\sigma_{\alpha}$ (up to a factor of $-2$). Now $\sigma_{\alpha}$ is a holomorphic two-form with respect to some complex structure, and therefore $\sigma_{\alpha}^{n+1}=0$ by degree reasons. Differentiation with respect to $\alpha$ yields $\sigma_{\alpha}^n\tau_{\alpha}=0$. Therefore, up to a constant, the $(n+k)$th term of $\Phi_{\alpha,\beta}$ is
\begin{eqnarray*}
(\alpha-\beta)^{-k}\sigma_{\alpha,\beta}^{n+k} & = & (\alpha-\beta)^{-k}(\sigma_{\alpha}+(\alpha-\beta)\tau_{\alpha})^{n+k} \\
 & = & (\alpha-\beta)^{-k}\sum_{l=k+1}^{n+k}{n+k\choose l}\sigma_{\alpha}^{n+k-l}(\alpha-\beta)^l\tau_{\alpha}^l \\
 & = & \sum_{l=k+1}^{n+k}{n+k\choose l}\sigma_{\alpha}^{n+k-l}(\alpha-\beta)^{l-k}\tau_{\alpha}^l,
\end{eqnarray*}
which is clearly divisible by $(\alpha-\beta)$.
\end{prf}

\begin{exm}
When $n=1$, we find that
\begin{eqnarray*}
(\sigma-(\alpha+\beta)\omega_I-\alpha\beta\bar{\sigma})^2 & = & (\alpha+\beta)^2\omega_I^2-2\alpha\beta\sigma\bar{\sigma} \\
 & = & (\alpha+\beta)^2\omega_I^2-4\alpha\beta\omega_I^2 \\
 & = & (\alpha-\beta)^2\omega_I^2
 \end{eqnarray*}
since $\sigma\bar{\sigma}=\omega_J^2+\omega_K^2=4\mathrm{vol}=2\omega_I^2$, while $\sigma^2$, $\sigma\omega_I$, $\omega_I\bar{\sigma}$, and $\bar{\sigma}^2$ all vanish for degree reasons. Substituting this in yields
$$\Phi_{\alpha,\beta}=\sigma-(\alpha+\beta)\omega_I+i(\alpha-\beta)(1-\mathrm{vol})-\alpha\beta\bar{\sigma}.$$
When $\alpha=\beta$ this gives the two-form $\sigma-2\alpha\omega_I-\alpha^2\bar{\sigma}$.
\end{exm}


\begin{rmk}
In addition to the diagonal $\Delta\subset\mathbb{CP}^1\times\mathbb{CP}^1$ parametrizing generalized complex structures $\mathcal{J}_{I_{\alpha}}$ of complex type, the graph of the antipodal map
$$\bar{\Delta}:=\{(\alpha,-\bar{\alpha}^{-1})\,|\,\alpha\in\mathbb{CP}^1\}\subset\mathbb{CP}^1\times\mathbb{CP}^1$$
parametrizes generalized complex structures $\mathcal{J}_{\omega_{\alpha}}$ of symplectic type. Indeed, $\bar{\Delta}$ is precisely the image under $f$ of the equators of the $\mathbb{CP}^1$-fibres of $\mathbb{P}(T\oplus\mathcal{O})$.
\end{rmk}

\subsection{Isotropic pure spinors}

Our second description of the $\mathbb{CP}^1\times\mathbb{CP}^1$-family of generalized complex structures only applies in real dimension four, i.e., when $n=1$. A hyperk{\"a}hler metric $g$ on a four-manifold $M$ determines a maximal positive subspace of $\mathrm{H}^2(M,\mathbb{R})$ with respect to the intersection pairing, namely the $3$-plane
$$V:=\langle\omega_I,\omega_J,\omega_K\rangle.$$
The holomorphic two-forms of the twistor family corresponding to $g$ lie in the complexification of $V$. Indeed, they are precisely the two-forms in $V\otimes\mathbb{C}$ that are isotropic with respect to the intersection pairing, since a $(2,0)$-form wedged with itself must vanish for degree reasons. Of course, these holomorphic two-forms are only defined up to a scalar, so the twistor family is really parametrized by
$$\left\{[X\omega_I+Y\omega_J+Z\omega_K]\in\mathbb{P}(V\otimes\mathbb{C})\,\left|\,(X\omega_I+Y\omega_J+Z\omega_K)^2=0\right.\right\}.$$
The equation reduces to $X^2+Y^2+Z^2=0$. We can solve by setting $X=-2\eta$, $Y=1-\eta^2$, and $Z=i(1+\eta^2)$, where $\eta$ is a parameter in $\mathbb{CP}^1$, and this recovers our earlier formula
$$\sigma_{\eta}=(\omega_J+i\omega_K)-2\eta\omega_I-\eta^2(\omega_J-i\omega_K).$$
In summary, the twistor family is parametrized by a conic in $\mathbb{P}(V\otimes\mathbb{C})\cong\mathbb{CP}^2$.

For generalized complex structures, our pure spinor can be a local section of $\wedge^{\mbox{even}}T^*\otimes\mathbb{C}$, rather than just $\wedge^2T^*\otimes\mathbb{C}$. (In fact, we could allow forms in all degrees, but we restrict to even degrees since we are interested in generalized complex structures of even type, such as complex and symplectic structures; recall that the parity of the type does not change in a family of generalized complex structures.) The hyperk{\"a}hler metric $g$ determines a maximal positive subspace of $\mathrm{H}^{\mbox{even}}(M,\mathbb{R})$ with respect to the Mukai pairing, namely the $4$-plane
$$W:=\left\langle\omega_I,\omega_J,\omega_K,1-\mathrm{vol}\right\rangle,$$
where $\mathrm{vol}$ denotes a volume form (equal to $\omega_I^2/2=\omega^2_J/2=\omega^2_K/2$). A spinor will be pure if it is isotropic with respect to the Mukai pairing, so
$$\left\{[X\omega_I+Y\omega_J+Z\omega_K+U(1-\mathrm{vol})]\in\mathbb{P}(W\otimes\mathbb{C})\,\left|\, (X\omega_I+Y\omega_J+Z\omega_K+U(1-\mathrm{vol}))^2=0\right.\right\}$$
gives a family of pure spinor lines. The equation reduces to $X^2+Y^2+Z^2+U^2=0$. We can solve by setting $X=-(\alpha+\beta)$, $Y=1-\alpha\beta$, $Z=i(1+\alpha\beta)$, and $U=i(\alpha-\beta)$, where $(\alpha,\beta)$ is a parameter in $\mathbb{CP}^1\times\mathbb{CP}^1$. In summary, we obtain a family of pure spinors
$$\Phi_{\alpha,\beta}=(\omega_J+i\omega_K)-(\alpha+\beta)\omega_I+i(\alpha-\beta)(1-\mathrm{vol})-\alpha\beta(\omega_J-i\omega_K)$$
parametrized by a quadric surface in $\mathbb{P}(W\otimes\mathbb{C})\cong\mathbb{CP}^3$. A priori, this will determine a $\mathbb{CP}^1\times\mathbb{CP}^1$-family of generalized almost complex structures on $M$. Moreover, $\omega_I$, $\omega_J$, $\omega_K$, and, of course, $1$ and $\mathrm{vol}$ are all $d$-closed. So for fixed $\alpha$ and $\beta$ we find $d\Phi_{\alpha,\beta}=0$, and Proposition~\ref{integrable} implies that these generalized almost complex structures are integrable. This $\mathbb{CP}^1\times\mathbb{CP}^1$-family of generalized complex structures on $M$ will be our generalized twistor family.

\begin{rmk}
The usual twistor family is parametrized by the hyperplane section given by intersecting the quadric surface with $\mathbb{P}(V\otimes\mathbb{C})\subset\mathbb{P}(W\otimes\mathbb{C})$. This of course gives a curve of bidegree $(1,1)$ in $\mathbb{CP}^1\times\mathbb{CP}^1$.
\end{rmk}

\begin{rmk}
The above approach is the one adopted by Huybrechts~\cite{huybrechts05}. He defines a generalized K3 surface as the underlying smooth four-manifold $M$ of a K3 surface equipped with a generalized K{\"a}hler structure (a pair of commuting generalized complex structures) together with choices of pure spinors defining this structure, normalized in a compatible way. Although we will mostly be interested in manifolds equipped with a single generalized complex structure, we shall see in the next subsection that $\mathbb{CP}^1\times\mathbb{CP}^1$ does indeed parametrize a family of generalized K{\"a}hler structures.
\end{rmk}

\subsection{Bi-Hermitian geometry}

Bi-Hermitian geometry was introduced by Gates, Hull, and Ro{\v c}ek~\cite{ghr84}. Its relation to generalized K{\"a}hler geometry was established by Gualtieri~\cite{gualtieri10}.

\begin{dfn}
An almost bi-Hermitian structure on a manifold $M$ consists of a metric $g$ together with a pair of almost complex structures, $I_+$ and $I_-$, that are compatible with the metric. We call $(g,I_+,I_-)$ a bi-Hermitian structure if both $I_+$ and $I_-$ are integrable.
\end{dfn}

\begin{dfn}
A generalized almost K{\"a}hler structure on a manifold $M$ consists of a pair of generalized almost complex structures, $\mathcal{J}$ and $\mathcal{J}^{\prime}$, which commute and such that the symmetric pairing
$$G(e_1,e_2):=-\langle\mathcal{J}e_1,\mathcal{J}^{\prime}e_2\rangle\qquad\mbox{for }e_1,e_2\in C^{\infty}(T\oplus T^*)$$
is positive definite, i.e., defines a generalized metric on $E=T\oplus T^*$. Note that $\mathcal{J}$ and $\mathcal{J}^{\prime}$ are automatically compatible with $G$, since they are compatible with $\langle\phantom{e},\phantom{e}\rangle$ and commute. We call $(\mathcal{J},\mathcal{J}^{\prime})$ a generalized K{\"a}hler structure if both $\mathcal{J}$ and $\mathcal{J}^{\prime}$ are integrable.
\end{dfn}

\begin{rmk}
One usually considers twisted structures, by adding a $B$-field to the bi-Hermitian structure and allowing Courant algebroids more general than $E=T\oplus T^*$. However, our application will only require untwisted structures.
\end{rmk}

\begin{exm}
If $(g,I,\omega)$ is a K{\"a}hler structure on $M$, then $(\mathcal{J}_I,\mathcal{J}_{\omega})$ is a generalized K{\"a}hler structure on $M$, with generalized metric $G$ as given at the start of this section.
\end{exm}

\begin{prp}
\label{bi-Hermitian}
Let $(g,I_+,I_-)$ be an almost bi-Hermitian structure on $M$. Then
$$\mathcal{J}=\frac{1}{2}\left(\begin{array}{cc} -(I_++I_-) & -(\omega_+^{-1}-\omega_-^{-1}) \\
  \omega_+-\omega_- & I_+^*+I_-^* \end{array}\right)$$
$$\mathcal{J}^{\prime}=\frac{1}{2}\left(\begin{array}{cc} -(I_+-I_-) & -(\omega_+^{-1}+\omega_-^{-1}) \\
  \omega_++\omega_- & I_+^*-I_-^* \end{array}\right)$$
give a generalized almost K{\"a}hler structure on $M$. If $(g,I_+,I_-)$ is a bi-Hermitian structure on $M$ and $g$ is K{\"a}hlerian with respect to both $I_+$ and $I_-$ then $(\mathcal{J},\mathcal{J}^{\prime})$ as defined above give a generalized K{\"a}hler structure on $M$.
\end{prp}

\begin{prf}
Gualtieri's Theorem~\cite{gualtieri10} is more general than this: he includes a $B$-field $b$ in the almost bi-Hermitian structure and determines the exact conditions on $(g,I_+,I_-,b)$ to ensure that $\mathcal{J}$ and $\mathcal{J}^{\prime}$, now also modified by a $B$-field transform, are integrable. When $b=0$ these conditions reduce to the integrability of $I_+$ and $I_-$, and the vanishing of $d\omega_+$ and $d\omega_-$, where $\omega_+$ and $\omega_-$ are the Hermitian forms of $(g,I_+)$ and $(g,I_-)$.
\end{prf}

\begin{exm}
A special case is when $I_+=I_-$. Then $(\mathcal{J},\mathcal{J}^{\prime})=(\mathcal{J}_{I_+},\mathcal{J}_{\omega_+})$ is the generalized K{\"a}hler structure coming from the K{\"a}hler structure $(g,I_+,\omega_+)$.
\end{exm}

On a hyperk{\"a}hler manifold $(M,g)$ we have a twistor family of complex structures, $I_{\eta}$ with $\eta\in\mathbb{CP}^1$. We can let both $I_+$ and $I_-$ vary in this family, setting $I_+=I_{\alpha}$ and $I_-=I_{\beta}$, leading to a family of almost bi-Hermitian structures $(g,I_{\alpha},I_{\beta})$ on $M$ parametrized by $(\alpha,\beta)\in\mathbb{CP}^1\times\mathbb{CP}^1$. By the above proposition, we get a family of generalized almost K{\"a}hler structures $(\mathcal{J}_{\alpha,\beta},\mathcal{J}^{\prime}_{\alpha,\beta})$ on $M$. Moreover, since $I_{\alpha}$ and $I_{\beta}$ are always integrable, and $g$ is K{\"a}hlerian with respect to $I_{\eta}$ for all $\eta\in\mathbb{CP}^1$, this is actually a family of generalized K{\"a}hler structures. The generalized metric $G$ is the one defined at the start of this section.

\begin{rmk}
Along the diagonal $\Delta\subset\mathbb{CP}^1\times\mathbb{CP}^1$,
$$(\mathcal{J}_{\alpha,\alpha},\mathcal{J}^{\prime}_{\alpha,\alpha})=(\mathcal{J}_{I_{\alpha}},\mathcal{J}_{\omega_{\alpha}})$$
comes from the K{\"a}hler structure $(I_{\alpha},\omega_{\alpha})$ on $M$. In other words, this subfamily of generalized K{\"a}hler structures is the usual twistor family, but equipped also with K{\"a}hler forms.
\end{rmk}

\begin{rmk}
Later we will show that $\mathcal{J}_{\alpha,\beta}$ depends holomorphically on $\alpha$ and $\beta$, in the sense that we can find a pure spinor defining $\mathcal{J}_{\alpha,\beta}$ which depends holomorphically on $\alpha$ and $\beta$. On the other hand, the pure spinor defining $\mathcal{J}^{\prime}_{\alpha,\beta}$ will depend on $\alpha$ and $\bar{\beta}$.
\end{rmk}

\subsection{Identifying the three families}

We now have three different descriptions of a family of generalized complex structures on a hyperk{\"a}hler manifold, each parametrized by $(\alpha,\beta)\in\mathbb{CP}^1\times\mathbb{CP}^1$. Of course, these are all the same family.

\begin{prp}
The three $\mathbb{CP}^1\times\mathbb{CP}^1$-families of generalized complex structures on $M$ described in the previous three subsections are actually the same family.
\end{prp}

\begin{prf}
We computed a pure spinor,
$$\Phi_{\alpha,\beta}=i^n(\alpha-\beta)^n\exp{\left(\frac{\sigma-(\alpha+\beta)\omega_I-\alpha\beta\bar{\sigma}}{i(\alpha-\beta)}\right)},$$
for the first family. When $n=1$ this agrees with the pure spinor of the second family (which only exists for $n=1$). It remains to identify the first and third families.

One easily observes that these familes agree along the diagonal, where the generalized complex structures are of complex type. Away from the diagonal, we saw that the first family is the $B$-field transform
$$e^{-B}\mathcal{J}_{\omega}e^B$$
of a generalized complex structure $\mathcal{J}_{\omega}$ coming from a symplectic structure. The real two-forms $B$ and $\omega$ can be expressed in terms of $\alpha$ and $\beta$. On the other hand, the third family is obtained by substituting
$$I_+=I_{\alpha}=\left(\frac{1-|\alpha|^2}{1+|\alpha|^2}\right)I+\left(\frac{2\re\alpha}{1+|\alpha|^2}\right)J+\left(\frac{2\im\alpha}{1+|\alpha|^2}\right)K,$$
$$I_-=I_{\beta}=\left(\frac{1-|\beta|^2}{1+|\beta|^2}\right)I+\left(\frac{2\re\beta}{1+|\beta|^2}\right)J+\left(\frac{2\im\beta}{1+|\beta|^2}\right)K,$$
$$\omega_+=\omega_{\alpha}=\left(\frac{1-|\alpha|^2}{1+|\alpha|^2}\right)\omega_I+\left(\frac{2\re\alpha}{1+|\alpha|^2}\right)\omega_J+\left(\frac{2\im\alpha}{1+|\alpha|^2}\right)\omega_K,$$
and
$$\omega_-=\omega_{\beta}=\left(\frac{1-|\beta|^2}{1+|\beta|^2}\right)\omega_I+\left(\frac{2\re\beta}{1+|\beta|^2}\right)\omega_J+\left(\frac{2\im\beta}{1+|\beta|^2}\right)\omega_K$$
into the first formula of Proposition~\ref{bi-Hermitian}. A lengthy calculation shows that the resulting generalized complex structure agrees with $e^{-B}\mathcal{J}_{\omega}e^B$; see Proposition~3.1.2 and Appendix~A of the first author's thesis~\cite{glover13}.
\end{prf}

It is now easy to find a pure spinor for the other family $\mathcal{J}^{\prime}_{\alpha,\beta}$ of generalized complex structures arising from the bi-Hermitian structure. We have used $\beta$ as a complex coordinate on $\mathbb{CP}^1$; let $\tilde{\beta}=\beta^{-1}$ be the coordinate in the other patch. We write $\mathcal{J}^{\prime}_{\alpha,\tilde{\beta}}$ for the same generalized complex structure $\mathcal{J}^{\prime}_{\alpha,\beta}$, expressed in coordinates $(\alpha,\tilde{\beta})$.

\begin{lem}
\label{Jprime}
Let
$$\Phi^{\prime}_{\alpha,\tilde{\beta}}:=i^n(\alpha+\bar{\tilde{\beta}})^n\exp{\left(\frac{\sigma-(\alpha-\bar{\tilde{\beta}})\omega_I+\alpha\bar{\tilde{\beta}}\bar{\sigma}}{i(\alpha+\bar{\tilde{\beta}})}\right)},$$
Then a pure spinor for $\mathcal{J}^{\prime}_{\alpha,\tilde{\beta}}$ is $\Phi^{\prime}_{\alpha,\tilde{\beta}}$. In particular, this generalized complex structure does not depend holomorphically on the parameters in $\mathbb{CP}^1\times\mathbb{CP}^1$. (It would depend holomorphically if we took the conjugate complex structure on the second factor, $\mathbb{CP}^1\times\overline{\mathbb{CP}^1}$.)
\end{lem}

\begin{prf}
Observe that replacing $I_-$ with $-I_-$ (and $\omega_-$ with $-\omega_-$) in Proposition~\ref{bi-Hermitian} interchanges $\mathcal{J}$ and $\mathcal{J}^{\prime}$. In our case $I_-=I_{\beta}$, so
$$-I_-=-I_{\beta}=I_{-\bar{\beta}^{-1}},$$
where $-\bar{\beta}^{-1}$ is the point on $\mathbb{CP}^1=S^2$ that is antipodal to $\beta$. Another way to state this is that $\mathcal{J}^{\prime}_{\alpha,\beta}$ is the pull back of $\mathcal{J}_{\alpha,\beta}$ by the map
\begin{eqnarray*}
\mathbb{CP}^1\times\mathbb{CP}^1 & \longrightarrow & \mathbb{CP}^1\times\mathbb{CP}^1 \\
(\alpha,\beta) & \longmapsto & (\alpha,-\bar{\beta}^{-1})
\end{eqnarray*}
given by the identity on the first factor and the antipodal map on the second. It follows that we obtain a pure spinor for $\mathcal{J}^{\prime}_{\alpha,\beta}$ by taking the pure spinor $\Phi_{\alpha,\beta}$ for $\mathcal{J}_{\alpha,\beta}$ and replacing $\beta$ by $-\bar{\beta}^{-1}$. Writing this in terms of the coordinate $\tilde{\beta}=\beta^{-1}$ then yields $\Phi^{\prime}_{\alpha,\tilde{\beta}}$.
\end{prf}

\subsection{Construction}

We have described a $\mathbb{CP}^1\times\mathbb{CP}^1$-family of generalized complex structures $\mathcal{J}_{\alpha,\beta}$ on $M$. Now we assemble them into a generalized complex structure on a single space.

\begin{thm}
Let $(M,g)$ be a hyperk{\"a}hler manifold. Denote by $\mathcal{X}$ the smooth manifold $M\times\mathbb{CP}^1\times\mathbb{CP}^1$ equipped with the generalized almost complex structure
$$\mathcal{J}:=\mathcal{J}_{\alpha,\beta}\times\mathcal{J}_{I_{\mathbb{CP}^1}}\times\mathcal{J}_{I_{\mathbb{CP}^1}},$$
where $I_{\mathbb{CP}^1}$ denotes the standard complex structure on $\mathbb{CP}^1$. Then $\mathcal{J}$ is integrable and thus $\mathcal{X}$ is a generalized complex manifold.
\end{thm}

\begin{dfn}
We call $\mathcal{X}$ the generalized twistor space of $(M,g)$.
\end{dfn}

\begin{prf}
A pure spinor for $\mathcal{J}_{\alpha,\beta}$ is $\Phi_{\alpha,\beta}$, while $d\alpha$ and $d\beta$ give pure spinors for $\mathcal{J}_{I_{\mathbb{CP}^1}}$ on the second and third factors. Combining these yields a pure spinor
$$\Psi:=\Phi_{\alpha,\beta}\wedge d\alpha\wedge d\beta$$
for $\mathcal{J}$. Now $\Phi_{\alpha,\beta}$ is holomorphic in $\alpha$ and $\beta$, and the other forms appearing in $\Phi_{\alpha,\beta}$, namely $\omega_I$, $\omega_J$, $\omega_K$, and their powers and products, are all $d$-closed. Therefore we find
$$d\Psi=(d\Phi_{\alpha,\beta})\wedge d\alpha\wedge d\beta=\frac{\partial\Phi_{\alpha,\beta}}{\partial\alpha}d\alpha\wedge d\alpha\wedge d\beta+
\frac{\partial\Phi_{\alpha,\beta}}{\partial\beta}d\beta\wedge d\alpha\wedge d\beta=0.$$
It follows from Proposition~\ref{integrable} that $\mathcal{J}$ is integrable.
\end{prf}

Since $\mathbb{CP}^1\times\mathbb{CP}^1$ parametrizes a family of generalized K{\"a}hler structures on $M$, we can construct a second generalized complex structure on $M\times\mathbb{CP}^1\times\mathbb{CP}^1$.

\begin{thm}
Define a generalized almost complex structure on $M\times\mathbb{CP}^1\times\mathbb{CP}^1$ by
$$\mathcal{J}^{\prime}:=\mathcal{J}^{\prime}_{\alpha,\beta}\times\mathcal{J}_{I_{\mathbb{CP}^1}}\times\mathcal{J}_{-I_{\mathbb{CP}^1}},$$
where $I_{\mathbb{CP}^1}$ denotes the standard complex structure on $\mathbb{CP}^1$ and $-I_{\mathbb{CP}^1}$ denotes the conjugate complex structure (the complex structure of $\overline{\mathbb{CP}}^1$). Then $\mathcal{J}^{\prime}$ is integrable. Together, $\mathcal{J}$ and $\mathcal{J}^{\prime}$ make the generalized twistor space $\mathcal{X}$ into a generalized pseudo-K{\"a}hler manifold, i.e., it satisfies the definition of a generalized K{\"a}hler except that the generalized metric $G$ is indefinite (and non-degenerate).
\end{thm}

\begin{prf}
We saw in Lemma~\ref{Jprime} that a pure spinor for $\mathcal{J}^{\prime}_{\alpha,\tilde{\beta}}$, in the coordinates $(\alpha,\tilde{\beta})$ on $\mathbb{CP}^1\times\mathbb{CP}^1$, is
$$\Phi^{\prime}_{\alpha,\tilde{\beta}}:=i^n(\alpha+\bar{\tilde{\beta}})^n\exp{\left(\frac{\sigma-(\alpha-\bar{\tilde{\beta}})\omega_I+\alpha\bar{\tilde{\beta}}\bar{\sigma}}{i(\alpha+\bar{\tilde{\beta}})}\right)}.$$
Pure spinors for $\mathcal{J}_{I_{\mathbb{CP}^1}}$ and $\mathcal{J}_{-I_{\mathbb{CP}^1}}$ on the second and third factors are given by $d\alpha$ and $d\bar{\tilde{\beta}}$, respectively. Therefore
$$\Psi^{\prime}:=\Phi^{\prime}_{\alpha,\tilde{\beta}}\wedge d\alpha\wedge d\bar{\tilde{\beta}}$$
is a pure spinor for $\mathcal{J}^{\prime}$. A calculation similar to the one in the previous proof shows that $d\Psi^{\prime}=0$, and then Proposition~\ref{integrable} implies that $\mathcal{J}^{\prime}$ is integrable.

It is easy to show that $\mathcal{J}$ and $\mathcal{J}^{\prime}$ commute. Recall that the generalized metric corresponding to the pair $(\mathcal{J},\mathcal{J}^{\prime})$ is given by
$$G(e_1,e_2):=\langle\mathcal{J}e_1,\mathcal{J}^{\prime}e_2\rangle.$$
This gives $g(X,Y)$ and $g^{-1}(\xi,\eta)$ for $X,Y\in TM$ and $\xi,\eta\in T^*M$, respectively, where $g$ is the hyperk{\"a}hler metric on $M$; so $G$ is positive definite on the first factor $M$. For $e_1,e_2\in T\mathbb{CP}^1\oplus T^*\mathbb{CP}^1$ for the second factor $\mathbb{CP}^1$, we find that
$$G(e_1,e_2)=\langle\mathcal{J}_{I_{\mathbb{CP}^1}}e_1,\mathcal{J}_{I_{\mathbb{CP}^1}}e_2\rangle=\langle e_1,e_2\rangle,$$
which has signature $(2,2)$. For $e_1,e_2\in T\mathbb{CP}^1\oplus T^*\mathbb{CP}^1$ for the third factor $\mathbb{CP}^1$, we find that
$$G(e_1,e_2)=\langle\mathcal{J}_{I_{\mathbb{CP}^1}}e_1,\mathcal{J}_{-I_{\mathbb{CP}^1}}e_2\rangle=-\langle e_1,e_2\rangle,$$
which again has signature $(2,2)$. Altogether, we see that $G$ is non-degenerate but not positive definite: it has signature $(8n+4,4)$.
\end{prf}

\begin{rmk}
Let us compare our construction to other constructions of generalized twistor spaces in the literature. Davidov and Muskarov~\cite{dm06} started with a $2n$-dimensional manifold $M$ and considered the bundle $\mathcal{G}\rightarrow M$ whose fibre over $p\in M$ parametrizes complex structures on $T_pM\oplus T^*_pM$. If $M$ is equipped with a linear connection, they described two generalized almost complex structures on the total space of $\mathcal{G}$, analogues of the Atiyah-Hitchin-Singer and Eells-Salamon almost complex structures. They proved that the first is integrable only when the dimension of $M$ is two or the linear conection is flat ($M$ is affine), while the second is never integrable. In a subsequent paper~\cite{dm07}, they showed that $\mathcal{G}$ admits a generalized K{\"a}hler structure if the dimension of $M$ is two {\em and\/} the linear connection is flat. The main difference to our construction is that their generalized twistor spaces are always non-compact, because the fibres $\mathcal{G}_p\cong\mathrm{SO}(2n,2n)/\mathrm{U}(n,n)$ parametrize complex structures compatible only with the indefinite inner product $\langle\phantom{e},\phantom{e}\rangle$. On the other hand, we introduced a Riemannian metric $g$ and the fibres of our generalized twistor space parametrize complex structures compatible with the associated generalized metric $G$, which is positive definite.

Pantilie~\cite{pantilie11} described a $\mathbb{CP}^1$-family of generalized complex structures on a holomorphic symplectic manifold $M$, and assembled them into a generalized complex structure on $M\times\mathbb{CP}^1$. This is the $\mathbb{CP}^1$-family that we described at the end of Section~2, and later used as the $\mathbb{CP}^1$-fibres of our $\mathbb{P}(T\oplus\mathcal{O})$-family of generalized complex structures on a hyperk{\"a}hler manifold. Therefore for hyperk{\"a}hler manifolds, Pantilie's twistor space is a submanifold of our generalized twistor space.

Finally, Bredthauer~\cite{bredthauer07} described generalized twistor spaces for hyperk{\"a}hler manifolds, using bi-Hermitian geometry as in our third approach. However, it is not clear how to derive a pure spinor via this approach: he states a formula for a pure spinor, but it seems to be incorrect. In short, we felt that this part of his paper, and consequently his proof of the integrability of the generalized almost complex structure, was incomplete.
\end{rmk}

\subsection{Properties}

We will describe some properties of the generalized twistor space $\mathcal{X}$, i.e., $M\times\mathbb{CP}^1\times\mathbb{CP}^1$ equipped with the generalized complex structure $\mathcal{J}$, that are analogues of the properties of the usual twistor space $Z$ described in Theorem~\ref{twistor}.

\begin{lem}
The generalized twistor space admits a fibration $p:\mathcal{X}\rightarrow\mathbb{CP}^1\times\mathbb{CP}^1$, in the sense that $\mathbb{CP}^1\times\mathbb{CP}^1$ is obtained from $\mathcal{X}$ via generalized reduction (for example, see Bursztyn, Cavalcanti, and Gualtieri~\cite{bcg06}).
\end{lem}

\begin{prf}
For the reduction of generalized complex structures described in~\cite{bcg06}, one starts with an extended action $\rho:\mathfrak{a}\rightarrow \Gamma(E_{\mathcal{X}})$ of a connected Lie group on the exact Courant algebroid $E_{\mathcal{X}}=T\mathcal{X}\oplus T^*\mathcal{X}$. In our case, we do not have a group, but defining $\mathfrak{a}:=E_M$ we have the inclusion
$$\rho:\mathfrak{a}\longrightarrow E_{\mathcal{X}}=E_M\oplus E_{\mathbb{CP}^1}\oplus E_{\mathbb{CP}^1},$$
and the `vertical' distribution $TM\subset T\mathcal{X}=TM\oplus T\mathbb{CP}^1\oplus T\mathbb{CP}^1$ gives a foliation $\mathcal{F}$ with quotient $\mathcal{X}/\mathcal{F}\cong\mathbb{CP}^1\oplus\mathbb{CP}^1$. These statements are all trivial, following automatically from the fact that $\mathcal{X}$ is a product as a smooth manifold. Less trivial, but true nonetheless, is the fact that the reduction procedure applied to $\mathcal{J}$ on $\mathcal{X}$ produces the generalized complex structure $\mathcal{J}_{I_{\mathbb{CP}^1}}\times\mathcal{J}_{I_{\mathbb{CP}^1}}$ associated to the standard complex structure on $\mathbb{CP}^1\times\mathbb{CP}^1$. This is a consequence of the block form of $\mathcal{J}=\mathcal{J}_{\alpha,\beta}\times\mathcal{J}_{I_{\mathbb{CP}^1}}\times\mathcal{J}_{I_{\mathbb{CP}^1}}$, which therefore preserves the decomposition of $E_{\mathcal{X}}=E_M\oplus E_{\mathbb{CP}^1}\oplus E_{\mathbb{CP}^1}$.
\end{prf}

\begin{rmk}
The moduli space of complex structures on $M$ is itself a complex manifold. The $\mathbb{CP}^1$ base of the twistor space is then a holomorphic curve inside this moduli space. In such a situation, the total space of the family will have a natural complex structure; this is the twistor space $Z$.

Similarly, Gualtieri~\cite{gualtieri11} showed that the moduli space of generalized complex structures on $M$ is itself a complex manifold. Our $\mathbb{CP}^1\times\mathbb{CP}^1$-family of generalized complex structures $\mathcal{J}_{\alpha,\beta}$ is holomorphic, in the sense that the pure spinor $\Phi_{\alpha,\beta}$ depends holomorphically on the parameters $\alpha$ and $\beta$. Therefore $\mathbb{CP}^1\times\mathbb{CP}^1$ is a complex submanifold of the moduli space of generalized complex structures, and again, the total space of this family should have a natural generalized complex structure; this is the generalized twistor space $\mathcal{X}$.
\end{rmk}

A submanifold $Y\subset X$ of a complex manifold is a complex submanifold if $TY\subset TX$ is preserved by the complex structure. This definition does not extend easily to generalized complex geometry, because $T^*Y$ is a quotient of $T^*X$, rather than a subbundle. Gualtieri defines the generalized tangent bundle $\tau_Y$ of a submanifold to be a certain extension of $TY$ by the conormal bundle $N^*Y$; if $\tau_Y\subset TX\oplus T^*X$ is preserved by the generalized complex structure then he calls $Y$ a generalized complex brane (see Section~6 of~\cite{gualtieri11} for details).

We wish to study the `holomorphic' sections of $p:\mathcal{X}\rightarrow\mathbb{CP}^1\times\mathbb{CP}^1$. These are {\em not\/} generalized complex branes. However, because of the block form of our generalized complex structure $\mathcal{J}$, we can define another kind of distinguished submanifold.

\begin{dfn}
Let $X$ be a generalized complex manifold given by a smooth product $M\times N$ and a generalized complex structure of block form, $\mathcal{J}=\mathcal{J}_1\times\mathcal{J}_2$, i.e., $\mathcal{J}_1$ is an endomorphism of $TM\oplus T^*M$ (which could depend on $n\in N$) and $\mathcal{J}_2$ is an endomorphism of $TN\oplus T^*N$ (which could depend on $m\in M$). For each $n\in N$, we call $M\times\{n\}\subset M\times N=X$ a complex factor manifold. Likewise, for each $m\in M$, we call $\{m\}\times N\subset M\times N=X$ a complex factor manifold.
\end{dfn}

\begin{rmk}
The point of this definition is that there is a canonical splitting
$$E_X=E_M\oplus E_N,$$
where $E_X:=TX\oplus T^*X$, etc., and $\mathcal{J}$ preserves the splitting because it is of block form. Thus $\mathcal{J}$ induces a generalized complex structure on each complex factor manifold.
\end{rmk}

\begin{exm}
The fibres of $p:\mathcal{X}\rightarrow\mathbb{CP}^1\times\mathbb{CP}^1$ are complex factor manifolds. Indeed, the fibre above $(\alpha,\beta)$ is precisely the hyperk{\"a}hler manifold $M$ equipped with the generalized complex structure $\mathcal{J}_{\alpha,\beta}$.
\end{exm}

\begin{lem}
For each $m\in M$, the section $Q_m:=\{m\}\times\mathbb{CP}^1\times\mathbb{CP}^1$ of $p:\mathcal{X}\rightarrow\mathbb{CP}^1\times\mathbb{CP}^1$ is a complex factor manifold. The `generalized normal bundle' $E_{\mathcal{X}}/E_{Q_m}$ of $Q_m$ is isomorphic to
$$\mathcal{O}(1,0)^{\oplus 2n}\oplus\mathcal{O}(0,1)^{\oplus 2n}.$$
\end{lem}

\begin{prf}
The fact that $Q_m$ is a complex factor manifold is immediate from the definition. The induced generalized complex structure is $\mathcal{J}_{I_{\mathbb{CP}^1}}\times\mathcal{J}_{I_{\mathbb{CP}^1}}$, where we identify $Q_m$ with $\mathbb{CP}^1\times\mathbb{CP}^1$. It remains to identify the generalized normal bundle.

As a smooth bundle, $E_{\mathcal{X}}/E_{Q_m}$ is the trivial bundle $(E_M)_m \times Q_m$. However, the varying complex structure $\mathcal{J}_{\alpha,\beta}|_{(E_M)_m}$ gives it a non-trivial holomorphic structure. We use the bi-Hermitian/generalized K{\"a}hler correspondence and adopt the notation of Section~2.3 of Gualtieri~\cite{gualtieri10}. Let
$$E_M\otimes\mathbb{C}=L_+\oplus\bar{L}_+$$
be the decomposition into $+i$ and $-i$-eigenbundles for $\mathcal{J}$. Then $E_M$ equipped with the complex structure $\mathcal{J}$ is isomorphic to $L_+$. We can further decompose
$$L_+=\ell_+\oplus\ell_-$$
into $+i$ and $-i$-eigenbundles for $\mathcal{J}^{\prime}$. As Gualtieri explains, there are isomorphisms
$$\ell_{\pm}\cong\{X\pm g(X)\,|\,X\in T^{1,0}_{\pm}M\},$$
where $T^{1,0}_{\pm}M$ is the $+i$-eigenbundle of $I_{\pm}$. For our generalized twistor space, $I_+=I_{\alpha}$ depends only on the first factor in $\mathbb{CP}^1\times\mathbb{CP}^1$. Moreover, for fixed $m\in M$, the holomorphic bundle $(T^{1,0}_+M)_m$ over this first factor $\mathbb{CP}^1$ is precisely the normal bundle of a holomorphic section of the usual twistor space $p:Z\rightarrow\mathbb{CP}^1$, which is isomorphic to $\mathcal{O}(1)^{\oplus 2n}$ by Theorem~\ref{twistor}(2). Because it is trivial over the second factor, we get $\mathcal{O}(1,0)^{\oplus 2n}$ over $Q_m=\mathbb{CP}^1\times\mathbb{CP}^1$. Similarly, the holomorphic bundle $(T^{1,0}_-M)_m$ over $Q_m=\mathbb{CP}^1\times\mathbb{CP}^1$ is isomorphic to $\mathcal{O}(0,1)^{\oplus 2n}$. Putting everything together, we find that
$$E_M\cong L_+=\ell_+\oplus\ell_-\cong T^{1,0}_+M\oplus T^{1,0}_-M,$$
and restricting to $Q_m$, the generalized normal bundle is
$$(E_M)_m\times Q_m\cong\mathcal{O}(1,0)^{\oplus 2n}\oplus\mathcal{O}(0,1)^{\oplus 2n}.$$
\end{prf}

\begin{lem}
There is a `holomorphic' section of $\Lambda^{\mbox{even}}T^*_F(n,n)$, where $T_F$ denotes the tangent bundle to the fibres, defining a pure spinor on each fibre of $p:\mathcal{X}\rightarrow\mathbb{CP}^1\times\mathbb{CP}^1$.
\end{lem}

\begin{prf}
This is just the family of pure spinors $\Phi_{\alpha,\beta}$. Recall that expanding the exponential gives
$$\Phi_{\alpha,\beta}=\sum_{j=0}^{2n}\frac{1}{j!}i^{n-j}(\alpha-\beta)^{n-j}(\sigma-(\alpha+\beta)\omega_I-\alpha\beta\bar{\sigma})^j.$$
It is clear that the $j\leq n$ terms have bidegree at most $(n,n)$ in $(\alpha,\beta)$. For $j>n$, the proof of Lemma~\ref{alpha=beta} shows that the $j=n+k$ term
$$(\alpha-\beta)^{-k}\sigma_{\alpha,\beta}^{n+k}=\sum_{l=k+1}^{n+k}{n+k\choose l}\sigma_{\alpha}^{n+k-l}(\alpha-\beta)^{l-k}\tau_{\alpha}^l$$
has degree at most $n$ in $\beta$. By symmetry, it must also have degree at most $n$ in $\alpha$. Therefore $\Phi_{\alpha,\beta}$ has bidegree $(n,n)$ in $(\alpha,\beta)$.
\end{prf}

\begin{lem}
The generalized twistor space $\mathcal{X}$ has a real structure compatible with the other structures and inducing the product of the antipodal maps, i.e.,
$$(\alpha,\beta)\longmapsto \left(-\bar{\alpha}^{-1},-\bar{\beta}^{-1}\right),$$
on $\mathbb{CP}^1\times\mathbb{CP}^1$.
\end{lem}

\begin{prf}
For the usual twistor space $Z=M\times\mathbb{CP}^1$ the real involution is given by
$$\tau:(m,\eta)\longmapsto\left(m,-\bar{\eta}^{-1}\right).$$
Compatibility with the holomorphic two-form follows from
$$\sigma_{-\bar{\eta}^{-1}}=-\bar{\eta}^{-2}\bar{\sigma}_{\eta}.$$
For the generalized twistor space $\mathcal{X}=M\times\mathbb{CP}^1\times\mathbb{CP}^1$ the real involution is given by
$$\tau:(m,\alpha,\beta)\longmapsto \left(m,-\bar{\alpha}^{-1},-\bar{\beta}^{-1}\right).$$
This map takes $I_+=I_{\alpha}$ to $I_{-\bar{\alpha}^{-1}}=-I_{\alpha}=-I_+$, and similarly, takes $I_-$ to $-I_-$. It therefore takes $\mathcal{J}$ to $-\mathcal{J}$, proving compatibility with the generalized complex structure ($\tau$ also takes $\mathcal{J}^{\prime}$ to $-\mathcal{J}^{\prime}$). The fibration $p:\mathcal{X}\rightarrow\mathbb{CP}^1\times\mathbb{CP}^1$ and holomorphic sections are clearly preserved by $\tau$. Finally, a calculation shows that
$$\Phi_{-\bar{\alpha}^{-1},-\bar{\beta}^{-1}}=\frac{(-1)^n}{\bar{\alpha}^n\bar{\beta}^n}\bar{\Phi}_{\alpha,\beta},$$
and therefore $\Phi_{\alpha,\beta}$ is also compatible with $\tau$.
\end{prf}

\begin{lem}
The hyperk{\"a}hler manifold $(M,g)$ can be recovered from its generalized twistor space $\mathcal{X}$, equipped with the structure described in the previous four lemmas, i.e., the fibration over $\mathbb{CP}^1\times\mathbb{CP}^1$, the holomorphic sections with generalized normal bundles as specified, the family of pure spinors, and the real structure.
\end{lem}

\begin{prf}
Let $p:\mathcal{X}\rightarrow\mathbb{CP}^1\times\mathbb{CP}^1$ be the fibration, let $\Delta$ be the diagonal in $\mathbb{CP}^1\times\mathbb{CP}^1$, and let $Z:=p^{-1}(\Delta)$. Then $p|_Z:Z\rightarrow\Delta\cong\mathbb{CP}^1$ is the usual twistor space of $(M,g)$. Moreover, the holomorphic sections of $p$ that are invariant under $\tau$ intersect $Z$ in rational curves with generalized normal bundles isomorphic to
$$\mathcal{O}(1,0)^{\oplus 2n}\oplus\mathcal{O}(0,1)^{\oplus 2n}|_{\Delta}\cong\mathcal{O}(1)^{\oplus 4n}.$$
(In the usual twistor space the normal bundles are isomorphic to $T^{1,0}_mM\cong\mathcal{O}(1)^{\oplus 2n}$, but then the generalized normal bundles will be isomorphic to $T^{1,0}_mM\oplus \Omega^{0,1}_mM\cong\mathcal{O}(1)^{4n}$.) These are the real twistor lines in $Z$. The family of pure spinors restricted to $Z$ look like $\sigma_{\eta}^n$, allowing us to recover the holomorphic section $\sigma_{\eta}$ of $\Lambda^2T^*_F(2)$ on $Z$. Finally, the real structure on $\mathcal{X}$ restricts to the real structure on $Z$. In short, the twistor space $Z$, equipped with its usual structure, sits inside the generalized twistor space $\mathcal{X}$. By Theorem~3.3 of Hitchin, et al.~\cite{hklr87}, we can recover $(M,g)$ from $Z$.
\end{prf}

\vspace*{3mm}
\begin{flushleft}
Department of Mathematics\hfill rglover3@z.rochester.edu\\
University of Rochester\hfill www.math.rochester.edu/people/faculty/rglover3\\
Rochester, NY 14627\\
USA\\
\vspace*{3mm}
Department of Mathematics\hfill sawon@email.unc.edu\\
University of North Carolina\hfill www.unc.edu/$\sim$sawon\\
Chapel Hill, NC 27599-3250\\
USA\\
\end{flushleft}

\end{document}